\documentclass[12pt]{elsarticle}

\usepackage{epic,eepic}
\usepackage{amssymb}
 \usepackage{amsthm}
\usepackage{mathrsfs}
\usepackage{fancybox}
\usepackage{fancyhdr}
\journal{Linear Algebra and its Applications}

\newtheorem{thm1}{Theorem}
\newtheorem{cor}[thm1]{Corollary}
\newtheorem{rem}[thm1]{Remark}
\newtheorem{deff}[thm1]{Definition}
\newtheorem{example}{Example}
\newtheorem{lem}[thm1]{Lemma}
\begin{document}

\begin{frontmatter}



\title{Matrices Totally Positive Relative to a Tree, II}


\author{R. S. Costas-Santos\corref{cor1}}
\ead{rscosa@gmail.com}
\ead[url]{http://www.rscosan.com}
\cortext[cor1]{Corresponding author}
\address{Dpto. de F\'isica y Matem\'aticas, 
Facultad de Ciencias, Universidad de Alcal\'a, 
28871 Alcal\'a de Henares, Spain}

\author{C. R. Johnson}
\ead{crjohnso@math.wm.edu}
\address{Department of Mathematics,
College of William and Mary, Williamsburg, VA 23187
}

\begin{abstract}
If $T$ is a labelled tree, a matrix $A$ is totally positive relative to
$T$, principal submatrices of $A$ associated with deletion of
pendent vertices of $T$ are P-matrices, and $A$ has positive
determinant, then the smallest absolute eigenvalue of $A$ is
positive with multiplicity 1 and its eigenvector is signed
according to $T$. This conclusion has been incorrectly 
conjectured under weaker hypotheses.
\end{abstract}

\begin{keyword}
Graph \sep Neumaier conjecture\sep Spectral theory \sep Sylvester's
identity \sep Totally positive matrix\sep Totally positive relative
to a tree.

\MSC[2010] Primary 05C50 \sep 15A18 \sep 15A48.
\end{keyword}

\end{frontmatter}


\section{Introduction} \label{intro-sec}
%%But before to prove such results we need to
%%introduce some preliminaries.

A real matrix is called \textit{totally positive}
(TP) if all its minors are positive, and it is a $P$-matrix 
if every principal minor is positive.

In \cite{JCB} the following weakening has been studied.  
An $n$-by-$n$ real matrix is {\it totally positive
relative to a given labelled tree $T$} on $n$ vertices (T-TP) if,
for each pair of pendent vertices $p$ and $q$ of
$T$, the matrix $A[\alpha]$ is TP when $\alpha$
is the ordered set of vertices of the unique induced
path of $T$ that connects $p$ and $q$.
%%Here $A[\alpha]$ is the principal
%%submatrix of $A$ with rows and columns $\alpha$,
%%re-ordered according to their labels.
If $T$ is a path with vertices labelled in order,
then TP and T-TP are the same. Note that we 
are going to refer to $T$ throughout as a labelled tree.

Of course, T-TP equivalently means that $A[\alpha]$
is TP for the vertices of any induced path of $T$,
as the unique path joining any pair of vertices of
$T$ is a subpath of some path joining pendent vertices.

It is known that a totally positive matrix has distinct
positive eigenvalues and that the smallest one has an
eigenvector that alternates in sign (see \cite{FJ} for
general background). Since a tree is bipartite,
there is a signing of the vertices so that neighbors have
different signs. For a labelled tree, $T$, let $\sigma$ be a
$\pm 1$ vector consistent with such a signing. We say that
$\sigma$ is {\it signed according to $T$}, and $\sigma$ is
unique up to multiplication by $\pm 1$. It had been
conjectured that if $A$ is T-TP, then  $A$ has a unique
absolute smallest real eigenvalue with an eigenvector
signed according to $T$. We call this the
{\it Neumaier conclusion}, after the original conjecture
by Arnold Neumaier, University of Vienna.
See \cite{JCB} for prior work.

This conjecture was proven for a few trees, but is false
in general. Here, our purpose
is to prove the original conjecture for all trees by adding
a hypothesis.
%%%%%%%%%%%%%%%%%%%%%%%%%%%%%%%%%%%%%%%%%%%%%%%%%%%%%%%%%%%%%%%
\section{Notation and Terminology} \label{not-term-sec}
Let us denote the set $\{1,\dots,n\}$ by $N$; Moreover, we 
will denote  by $N_{i}$ (resp. $N_{i,j}$, and
$N_{i,j,k}$) the set $N\setminus\{i\}$ (resp.
$N\setminus\{i,j\}$, and $N\setminus\{i,j,k\}$).

Let $A\in M_n(\mathbb R)$. For any
ordered index lists $\alpha$, $\beta \subseteq N$,
with $|\alpha|=|\beta|=k$, by $A[\alpha ; \beta]$ we mean
the $k$-by-$k$ submatrix of $A$ that lies in the
rows indexed by $\alpha$ and
the columns indexed by $\beta$, and with the
order of the rows (resp. columns) determined
by the order in $\alpha$ (resp. $\beta$);
by $A[\alpha]$ we mean $A[\alpha;\alpha]$;
by $A(i;j)$ we mean the $(n-1)$-by-$(n-1)$ submatrix
of $A$ that lies in the rows 
indexed by $N_i$ and  the columns indexed
by $N_j$;  and by $A(i)$ we mean
$A(i;i)$.

Suppose that $T$ is a labelled tree on
$n$ vertices.
If $\mathscr P$ is an induced path of $T$,
by $A[\mathscr P]$ we mean $A[\alpha]$ in which
$\alpha$ consists of the indices of the vertices
of $\mathscr P$ in the order in which they appear
along $\mathscr P$.
Since everything we discuss is independent of
reversal of order,  there is no ambiguity regarding intended 
direction.
\begin{deff} \label{def2.1}
For a given labelled tree $T$ on $n$
vertices, we say that $A \in M_{n}(\mathbb{R})$ is
T-TP if $A[\mathscr P]$ is TP
for each path  $\mathscr P$ connecting any
two pendent vertices.
\end{deff}

%%If a path is labelled in some other way,
%%a T-TP matrix is permutation similar to a TP
%%matrix.

Observe that for a T-TP matrix, properly less
is required than for a TP matrix; however,
like TP matrices, T-TP matrices are entry-wise
positive.

%%If $T$ is a tree but not a path, only
%%certain (re-ordered) proper principal
%%submatrices of a T-TP matrix are required to
%%be TP and the matrix itself need not be
%%permutation similar to a TP matrix.
\begin{deff} \label{def2.2}
For a given labelled tree $T$ on $n$ vertices, we say
that $A \in M_{n}(\mathbb{R})$ is {\tt pendent-$P$
relative to $T$} if all principal submatrices, associated
with the deletion of pendent vertices, are P-matrices.
\end{deff}

Note that since in a P-matrix all the principal minors are
positive the  property of being pendent-$P$
relative to a tree is preserved by permutation similarity.

\begin{deff} \label{def2.3}
For a given labelled tree $T$ on $n$ vertices, we say
that $A \in M_{n}(\mathbb{R})$ is {\tt $T$-positive}
if it is T-TP and pendent-$P$ relative to $T$.
\end{deff}

Our arguments strongly use the adjoint
of a T-TP matrix (or one satisfying additional
hypotheses), as a surrogate for the inverse, and
we frequently use Sylvester's determinantal
identity, along with ad hoc arguments,
to determine the sign pattern of the adjoint.

The version of {\it Sylvester's identity} we shall
use is the following \cite[(0.8.6.1)]{HoJo13}:
\begin{equation} \label{2.1}
\det A[\alpha;\beta]=\frac{\det A[\alpha';\beta']
\det A['\alpha;'\beta]-\det A[\alpha';'\beta]
\det A['\alpha;\beta']}{\det A['\alpha';'\beta']},
\end{equation}
in which $\alpha$ and $\beta$ are index sets of the
same size, $\alpha'$ (resp. $\beta'$) is $\alpha$
(resp. $\beta$) without the last index;
$'\!\alpha$ (resp. $'\!\beta$) is $\alpha$
(resp. $\beta$) without the first index, and
$'\!\alpha'$ (resp. $'\!\beta'$) is $\alpha$ (resp.
$\beta$) without the first index and last index.
%%We also adopt the notation that a \verb+^+ over an
%%index in an index set means that the index is
%%omitted from the set.
Note that, above, as throughout, these index sets
are ordered.
We also denote by $\widetilde A=(\widetilde a_{ij})$ the
adjoint of $A$.
%%%%%%%%%%%%%%%%%%%%%%%%%%%%%%%%%%%%%%%%%%%%%%%%%%%%%%%%%%%%%%%
\section{Main Result} \label{main-sec}
Our purpose here is to give hypotheses sufficient to
achieve the Neumaier conclusion relative to any tree.
Our approach is to give hypotheses so that $SA^{-1}S$
is an entry-wise positive matrix when $S$ is the signature matrix
determined by $\sigma$ signed according to $T$.
By Perron's Theorem this means that the smallest
eigenvalue is positive and has an eigenvector signed
according to $T$.
To this end our first result is.

\begin{thm1}\label{theo1}
Let $T$ be a labelled tree on $n$ vertices and $A\in M_n(\mathbb R)$
be T-positive. Then, if $\det A>0$, we have
$$
{\rm sign}(\det A(i;j))=(-1)^{i+j}\sigma_i \sigma_j
$$
in which $\sigma$ is signed according to $T$.
\end{thm1}
\begin{rem}
It is important to point out the fact that
$(-1)^{i+j}\det A(j;i)$ is the $(i,j)$ entry in
the adjoint matrix of $A$, i.e.,
$$
\det A(j;i)=(-1)^{i+j}\,\widetilde a_{ij}.
$$
For this reason we will write
$\widetilde a_{ij}$ instead of of $(-1)^{i+j}\det A(j;i)$
throughout the paper.
\end{rem}

Now, let $S_\sigma={\rm diag}(\sigma_1,\sigma_2,
\dots, \sigma_n)$ with $\sigma$ signed according
to $T$. We have

\begin{cor} \label{cor3.3}
If $T$ is a tree on $n$ vertices and $A\in
M_n(\mathbb R)$ is T-positive, with $\det A>0$, then
$$
S_\sigma A^{-1}  S_\sigma \ \ \mbox{is entry-wise positive}.
$$
Therefore, $A$ satisfies the Neumaier conclusion.
\end{cor}
%%%%%%%%%%%%%%%%%%%%%%%%%%%%%%%%%%%%%%%%%%%%%%%%%%%%%%%%%%%%%%%
\section{Supporting Facts and Proofs} \label{facts-sec}
In this section we give the results
that we need in order to prove Theorem
\ref{theo1}. We also deduce the corollaries from it.
First we state a technical result we
need to prove Lemma \ref{lemma4.4}.

\begin{lem} \label{lemma4.1}
Given a matrix $A\in M_n(\mathbb R)$, then for any three
distinct integers $i$, $j$, $k$, with $1\le i,j,k\le n$,
we  have
$$
\widetilde a_{ki}\! \det A[i,N_{i,j,k};\!i,N_{i,j,k}]+
\widetilde a_{kj}\! \det A[j,N_{i,j,k};\!i,N_{i,j,k}]+
\widetilde a_{kk}\! \det A[k,N_{i,j,k};\!i,N_{i,j,k}]\!=\!0.
$$
\end{lem}
\begin{proof}
Without loss of generality we can assume that we can 
assume that $1\le i<j<k\le n$.
To simplify the expressions we are going to denote
$N_{i,j,k}$ by $\alpha$. By Sylvester's identity (\ref{2.1}) 
and taking into account how the sign of the determinant 
changes after a permutation of the indices we get
$$\begin{array}{rl}
\widetilde a_{kj}=&(-1)^{k+j}(-1)^{2i-2}\det A[i,N_{i,j};
i,N_{i,k}]\!=\!(-1)^{n-1}\det A[i,N_{i,j,k},k;j,i,
N_{i,j,k}]\\[3mm]=& (-1)^n \displaystyle\frac
{\det A[i,\alpha;i,\alpha]\det A[\alpha,k;j,i,\alpha']
-\det A[i,\alpha;j,i,\alpha']\det A[\alpha,k;i,\alpha]}
{\det A[\alpha;i,\alpha']}\\[3mm] =& \displaystyle
\frac{\det A[i,\alpha;i,\alpha]\det A[\alpha,k;j,i,
\alpha']-\det A[i,\alpha;j,i,\alpha']\det A[\alpha,k;i,
\alpha]}{\det A[r,\alpha';i,\alpha']},
\end{array}$$
where $r$ is the last entry of $\alpha$.

On the other hand and following the same ideas
as before, we get
$$\begin{array}{rl}
\widetilde a_{ki}=&(-1)^{k+i+1}\det A[j,N_{i,j};
j,N_{j,k}]=(-1)^n \det A[j,\alpha,k;j,i,\alpha]
\\[3mm]=& \displaystyle(-1)^n\frac{\det A[j,\alpha;
j,i,\alpha']\det A[\alpha,k;i,\alpha]-\det A[j,\alpha;i,
\alpha]\det A[\alpha,k;j,i,\alpha']}{\det A[\alpha;i,
\alpha']}\\[3mm]=& \displaystyle\frac{\det A[j,\alpha;
j,i,\alpha']\det A[\alpha,k;i,\alpha]-\det A[j,\alpha;i,\alpha]
\det A[\alpha,k;j,i,\alpha']}{\det A[r,\alpha';i,\alpha']}.
\end{array}$$
Thus using the last two expressions and combining them
properly, we get
$$\begin{array}{rl}
\widetilde a_{ki}\det A[i,\alpha;i,\alpha]=& \displaystyle
-\det A[j,\alpha;i,\alpha]\left(\widetilde a_{kj}+
\frac{\det A[i,\alpha;j,i,\alpha']\det A[\alpha,k;i,\alpha]}
{\det A[r,\alpha';i,\alpha']}\right)\\[3mm] &
\displaystyle +\frac{\det A[j,\alpha;j,i,\alpha']
\det A[\alpha,k;i,\alpha]\det A[i,\alpha;i,\alpha]}
{\det A[r,\alpha';i,\alpha']}\\[3mm] =&  \displaystyle
-\widetilde a_{kj}\det A[j,\alpha;i,\alpha]-\frac{\det A[\alpha,
k;i,\alpha]}{\det A[r,\alpha';i,\alpha']}\big(\det A[i,\alpha;
j,i,\alpha']\\[3mm] & \times \displaystyle
\det A[j,\alpha;i,\alpha]\!-\!\det A[j,\alpha;j,i,
\alpha']\det A[i,\alpha;i,\alpha]\big)\\[3mm]
=& -\widetilde a_{kj}\det A[j,\alpha;i,\alpha]-\widetilde a_{kk}
\det A[k,\alpha;i,\alpha].
\end{array} $$
\end{proof}
It is important to point out that, via permutation
similarity, the labelling of the tree, per se, is
not important.
If the conjecture were correct for one labelling
of a given tree, it would be correct for another.
Indeed, it is an easy exercise to see that if a
path is labelled in some other way than consecutively,
a T-TP matrix still has the ``last'' eigenvector
signed according to the alternatively labelled path.

Once the next three lemmata are proven, Theorem
\ref{theo1} follows. In the first lemma we prove
the statement of the theorem for any two pendent
vertices. If the tree is not a path, then it has
at least 3 pendent vertices.
Then we prove the statement of the theorem
assuming $i$ is pendent and $j$ is any vertex,
and  in the last lemma we prove the the statement
of the theorem without assuming $i$ and $j$ are
pendent vertices.

We are going to prove these lemmata by induction
on the number of vertices $n$, $n\ge 2$,
of the tree $T$.
The cases $2\le n\le 4$ were proven in \cite{JCB},
so we will assume $n\ge 4$ and that $T$ is not a
path (in which case the claim is inmediate).
Recall that $\sigma$ is signed according to $T$. Then,
we need to prove
\begin{equation} \label{3.1}
{\rm sign}(\det A(i;j))=(-1)^{i+j}\sigma_i \sigma_j,
\end{equation}
for all $1\le i,j\le n$.
Note that if (\ref{3.1}) holds, since
$$
{\rm sign}(\det A(i;j))=(-1)^{i+j}\sigma_i \sigma_j
\ \iff \ {\rm sign}(\widetilde a_{ij})=\sigma_i\sigma_j
$$
the matrix
$$
{\rm diag}(\sigma_1,\cdots, \sigma_N)\,\widetilde A \, {\rm diag}
(\sigma_1,\cdots, \sigma_N)
$$
is entry-wise positive.

\begin{lem} \label{lemma4.2}
Under the same assumptions as in Theorem \ref{theo1}, for
any two different pendent vertices, $p_1$ and $p_2$,
$$
{\rm sign}(\det A(p_1;p_2))=(-1)^{p_1+p_2}\sigma_{p_1}
\sigma_{p_2}.
$$
\end{lem}
\begin{proof}
Since after removing a pendent vertex
of a tree it is still a tree (the tree has at
most 4 vertices and it is not a path), we can apply
the induction hypothesis to obtain
$$
\det A(p_1;p_2)=
\det A[N_{p_1},N_{p_2}]=(-1)^{p_1+p_2-1}\det A[p_2,N_{p_1,p_2};
p_1,N_{p_1,p_2}].
$$
Without loss of generality, let $p_3$ be  the last pendent vertex 
in $N$, with $p_3>\max\{p_1,p_2\}$. Therefore, if we denote
$N_{p_1,p_2,p_3}\cup\{p_3\}$ by $\alpha$ and use Sylvester's
identity we get that $(-1)^{p_1+p_2}\det A(p_1;p_2)$ is equal
to
$$
\frac{\det A[p_2,\alpha';\alpha]\det
A[\alpha;p_1,\alpha']-\det A[p_2,\alpha';p_1,\alpha']
\det A[\alpha;\alpha]}{\det A[\alpha';\alpha']}.
$$
%\begin{rem}
Notice that, since the tree has at least 3 pendent
vertices, we have rearranged the entries of $\alpha$
in such a way that the last element of $\alpha$ is the
pendent vertex $p_3$, i.e. $\alpha'\cup \{p_3\}=\alpha$;
while, for example, $\widetilde a_{p_1,p_2}|_{p_3}$
represents the entry $(p_1,p_2)$ of the adjoint
of the $(n-1)\times (n-1)$ submatrix of $A$ from which
the $p_3$-th row and $p_3$-th column are removed. By the
induction hypothesis sign$(\det A(p_3)(p_2;p_1))
=(-1)^{p_1+p_2}\sigma_{p_1}\sigma_{p_2}$.
%\end{rem}

Here the denominator is positive because $A$ is
pendent-$P$ and $p_3$ is a pendent vertex; the numerator
has the desired sign since (let us assume, for example, that $p_1<p_2$)
\begin{eqnarray*}
\mbox{sign}(\det A[p_2,\alpha';\alpha])
=&(-1)^{p_3}\sigma_{p_2}\sigma_{p_3},\\
\mbox{sign}(\det A[\alpha;p_1,\alpha'])
=&(-1)^{p_3}\sigma_{p_1}\sigma_{p_3}, \\
\mbox{sign}(\det A[p_2,\alpha';p_1,\alpha'])
=&-\sigma_{p_1}\sigma_{p_2}, \\
\mbox{sign}(\det A[\alpha;\alpha])
=&+.
\end{eqnarray*}

Observe that  if $p_1<p_2$ and due to the re-labeling after the deleting of 
$p_2$ in the new tree there is a shift in the resulting sign of 
$\det A(p_2)(p_1;p_3)$.

Then, since $p_1$, $p_2$, and $p_3$  are pendent vertices,
again by the induction hypothesis, we have
$$
{\rm sign}(\det A[p_2,\alpha';\alpha]
\det A[\alpha;p_1,\alpha'])=\sigma_{p_2}\sigma_{p_3}
\sigma_{p_1}\sigma_{p_3}=\sigma_{p_2}\sigma_{p_1},
$$
and
$$
-{\rm sign}(\det A[p_2,\alpha';p_1,\alpha'])
=\sigma_{p_1}\sigma_{p_2},
$$
so that the claim follows.
\end{proof}

Next, by using Lemma \ref{lemma4.1}, we are going to
prove the following result:
\begin{lem} \label{lemma4.4}
Under the same assumptions as in theorem \ref{theo1}, for
any pendent vertex $p$ and for any $i$, $1\le i\le N$,
$$
{\rm sign}(\det A(i;p))=(-1)^{i+p}\sigma_i \sigma_{p}.
$$
\end{lem}
\begin{proof}
If $i$ is a pendent vertex, $i\ne p$, the result
follows from Lemma \ref{lemma4.2}. If $i=p$ the
result follows since $p$ is a pendent vertex. Therefore,
by the pendent-$P$ hypothesis relative to $T$, we have
$$
\det A(p;p)=\det A(p)>0, \ \  \sigma_p\, \sigma_p>0 \qquad
\Rightarrow
\quad {\rm sign}(\det A(p;p))=\sigma_p \sigma_{p}.
$$
On the other hand, if $i$ is not a pendent vertex,
then setting in Lemma \ref{lemma4.1} the vertex $j$ as another
pendent vertex, namely $q$, and $k=p$, we get
$$
0\!=\!\widetilde a_{pi}\det \!A[i,N_{i,q,p};i,N_{i,q,p}]
+\widetilde a_{pq}\det \!A[q,N_{i,q,p};i,N_{i,q,p}]
+\widetilde a_{pp}\det \!A[p,N_{i,q,p};i,N_{i,q,p}].
$$
Taking into account that $p$ is a pendent vertex, by
hypothesis and induction, we have
$$\begin{array}{rl}
\det A[i,N_{i,q,p};i,N_{i,q,p}]=&\det A(p)(q;q)
>0, \\[3mm]
{\rm sign}(\det A[q,N_{i,q,p};i,N_{i,q,p}])=&-\sigma_q\,\sigma_i,\\[3mm]
{\rm sign}(\det A[p,N_{i,q,p};i,N_{i,q,p}])=&-\sigma_p\,\sigma_i,
\end{array}$$
$\widetilde a_{pp}>0$, and ${\rm sign}(\widetilde a_{pq})
=\sigma_p\,\sigma_q$.
Therefore
$$
{\rm sign}(\widetilde a_{pq}\det A[q,N_{i,q,p};i,N_{i,q,p}])
=-\sigma_p\,\sigma_i={\rm sign} \big(\widetilde
a_{p,p}\det A[p,N_{i,q,p};i,N_{i,q,p}]\big),
$$
and since $k=p$ and taking into account the equality of
Lemma \ref{lemma4.1}, we have
$$
(-1)^{k+i} {\rm sign} \det A(i;k)={\rm sign} \big(\widetilde
a_{k,k}\det A[k,N_{i,q,k};i,N_{i,q,k}]\big)=\sigma_k \sigma_i,
$$
so that the claim follows.
\end{proof}

For the last lemma we need to use {\it Jacobi's
identity} \cite[(0.8.4.1)]{HoJo13}
\begin{equation} \label{4.1}
\det A[\alpha;\beta]=(-1)^{p(\alpha,\beta)} \det A \
\det A^{-1}[N\setminus \beta;N\setminus \alpha],
\end{equation}
in which $|\alpha|=|\beta|$, and $p(\alpha,\beta)=
\sum_{i\in \alpha} i+\sum_{j\in \beta} j$.

\begin{lem} \label{lemma4.5}
Under the same assumptions as in Theorem \ref{theo1}, for
any pair $(i, j)$, neither of which is pendent,
$$
{\rm sign}(\det A(i;j))=(-1)^{i+j}\sigma_i \sigma_{j}.
$$
\end{lem}
\begin{proof}
We prove this by contradiction.
If we assume that ${\rm sign}(\det A(i;j))\ne (-1)^{i+j}
\sigma_i \sigma_{j}$ then, ${\rm sign}(\widetilde a_{ji})\ne
\sigma_i \sigma_{j}$. Let $p$ any pendant vertex. 
Then, on one hand, we
have
$$
\det \widetilde A[j,p;i,p]=\left|
\begin{array}{cc} \widetilde a_{ji} &
\widetilde a_{jp} \\ \widetilde a_{p,i} &
\widetilde a_{p,p}\end{array}\right|=
\widetilde a_{ji}\widetilde a_{p,p}-
\widetilde a_{jp}\widetilde a_{p,i},
$$
so that, by lemmata \ref{lemma4.2} and
\ref{lemma4.4}, we get
$$
{\rm sign} \left(\det \widetilde A[j,p;i,p]\right)
=-\sigma_i\sigma_j.
$$
On the other hand, since $\det A>0$, applying
{\it Jacoibi's identity} we have
$$\begin{array}{rl}
{\rm sign}\big(\det \widetilde A[j,p;i,p]\big)& ={\rm
sign}\big((-1)^{i+j}\det A[N_{i,p};N_{j,p}]\big)\\ &
={\rm sign}\big((-1)^{i+j}\det A(p)(i;j)\big),
\end{array}$$
so that it is equal to, by the induction hypothesis,
$\sigma_i\sigma_j$ which is a contradiction. Hence
the result follows.
\end{proof}

Now, because of the relationship between the
$(n-1)$-by-$(n-1)$ minors of $A$ and $\widetilde A$,
Theorem \ref{theo1} follows from lemmata \ref{lemma4.2},
\ref{lemma4.4}, and \ref{lemma4.5}, as all types of minors
are covered.
As $\det(A)>0$, because of the relation between $A^{-1}$
and $\widetilde A$, corollary \ref{cor3.3} follows.
Then, as the Perron root of $A^{-1}$ is the reciprocal of
the smallest absolute eigenvalue of $A$, that smallest
eigenvalue is positive and has multiplicity 1. Because of the
effect of similarity on eigenvectors (see \cite{HoJo13}) 
the result about the signing of its
eigenvector follows.
%%%%%%%%%%%%%%%%%%%%%%%%%%%%%%%%%%%%%%%%%%%%%%%%%%%%%%%%%%%%%%%
\section{Remarks} \label{outlook-sec}
We have shown that certain conditions on a matrix $A$, relative
to a tree, are sufficient to reach the Neumaier conclusion.
These conditions are more, see \cite{JCB}, than originally
conjectured, but the originally conjectured conditions (T-TP)
were not sufficient in general. We do not know if some of the
additional hypotheses can be omitted. It is difficult to
construct appropriate examples.

However, we do have some informative examples. It is possible
for matrix $A$ to be T-positive but have negative determinant and
 satisfy the Neumaier conclusion. We still do not know
how common this is.
\begin{example} \label{ex1}
For this example we have considered the 5-star and the following
5-by-5 matrix. It is easy to check that $A$ is pendent-$P$
relative to this tree and $\det (A)<0$.
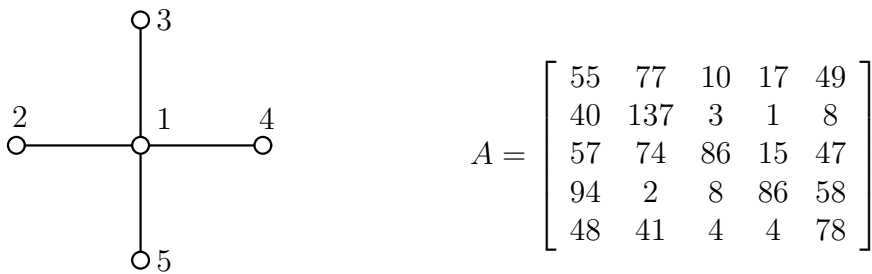
\begin{center}
\begin{minipage}{60mm}
\setlength{\unitlength}{0.00083333in}
\begingroup\makeatletter\ifx\SetFigFont\undefined%
\gdef\SetFigFont#1#2#3#4#5{%
  \reset@font\fontsize{#1}{#2pt}%
  \fontfamily{#3}\fontseries{#4}\fontshape{#5}%
  \selectfont}%
\fi\endgroup%
{\renewcommand{\dashlinestretch}{30}
\begin{picture}(1658,1727)(0,-10)
\thicklines
\put(53,812){\circle{106}}
\put(828,812){\circle{106}}
\put(1590,812){\circle{106}}
\put(828,95){\circle{106}}
\put(828,1595){\circle{106}}
\drawline(116,812)(770,812)%(764,812)
\drawline(881,812)(1539,812)%(1539,812)
\drawline(828,867)(828,1535)%(828,1535)
\drawline(828,150)(828,753)%(816,773)
\put(28,925){\makebox(0,0)[lb]{{\SetFigFont{12}{14.4}{\rmdefault}{\mddefault}{\updefault}2}}}
\put(928,25){\makebox(0,0)[lb]{{\SetFigFont{12}{14.4}{\rmdefault}{\mddefault}{\updefault}5}}}
\put(1566,912){\makebox(0,0)[lb]{{\SetFigFont{12}{14.4}{\rmdefault}{\mddefault}{\updefault}4}}}
\put(928,912){\makebox(0,0)[lb]{{\SetFigFont{12}{14.4}{\rmdefault}{\mddefault}{\updefault}1}}}
\put(928,1532){\makebox(0,0)[lb]{{\SetFigFont{12}{14.4}{\rmdefault}{\mddefault}{\updefault}3}}}
\end{picture}
}
\end{minipage}
%%%%%%%
\begin{minipage}{56mm}
\begin{eqnarray*}
A
=
\left[
\begin{array}{ccccc}
 55 & 77 & 10 & 17 & 49 \\
 40 & 137 & 3 & 1 & 8 \\
 57 & 74 & 86 & 15 & 47 \\
 94 & 2 & 8 & 86 & 58 \\
 48 & 41 & 4 & 4 & 78
\end{array}
\right]
\end{eqnarray*}
\end{minipage}
\end{center}
Note that in this example, the eigenvector associated
with the smallest eigenvalue, $\lambda_5 \sim -0.23$,
has the predicted sign pattern. Here is the eigenvector in
question, with each entry approximated to the nearest
hundredth:
\begin{eqnarray*}
\textbf{x}
\approx
\left[
  \begin{array}{c}
    -2.3 \\
    0.6 \\
    0.15 \\
    1.8 \\
    1 \\
  \end{array}
\right].
\end{eqnarray*}
The adjoint of $A$ is
$$
\widetilde A=\left[
\begin{array}{ccccc}
 70451860 & -27857784 & -4763560 & -11372966 & -30073840 \\
 -18274672 & 7046528 & 1241168 & 2950496 & 7815680 \\
 -4532012 & 1908264 & 18096 & 774494 & 2064504 \\
 -55473260 & 21866360 & 3770144 & 8668470 & 23888344 \\
 -30671880 & 12220096 & 2084744 & 4963592 & 12765448
\end{array}
\right].
$$
Both $\textbf{x}$ and $\widetilde A$  have the predicted sign pattern.
\end{example}

However, if $A$ is T-TP but not pendent-$P$ relative to T, the
Neumaier conclusion may fail.
\begin{example} \label{ex2} \cite{JCB} For this example we have considered
the following tree with 5 vertices and the following
5-by-5 matrix. It is easy to check that $\det A(5)<0$ therefore
$A$ is not pendent-$P$ relative to this tree, and $\det (A)<0$.
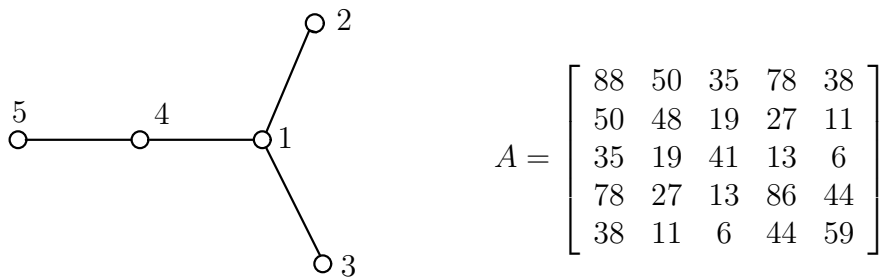
\begin{center}
\begin{minipage}{60mm}
\setlength{\unitlength}{0.00083333in}
\begingroup\makeatletter\ifx\SetFigFont\undefined%
\gdef\SetFigFont#1#2#3#4#5{%
  \reset@font\fontsize{#1}{#2pt}%
  \fontfamily{#3}\fontseries{#4}\fontshape{#5}%
  \selectfont}%
\fi\endgroup%
{\renewcommand{\dashlinestretch}{30}
\begin{picture}(2142,1691)(0,-10)
\thicklines
\put(68,850){\circle{106}}
\put(828,850){\circle{106}}
\put(1590,850){\circle{106}}
\put(1969,75){\circle{106}}
\put(1918,1578){\circle{111}}
\drawline(116,850)(774,850)
\drawline(891,850)(1539,850)
\drawline(1616,897)(1884,1532)(1928,1521)
\drawline(1616,798)(1952,126)
\put(28,963){\makebox(0,0)[lb]{{\SetFigFont{12}{14.4}
{\rmdefault}{\mddefault}{\updefault}5}}}
\put(1686,800){\makebox(0,0)[lb]{{\SetFigFont{12}{14.4}
{\rmdefault}{\mddefault}{\updefault}1}}}
\put(916,950){\makebox(0,0)[lb]{{\SetFigFont{12}{14.4}
{\rmdefault}{\mddefault}{\updefault}4}}}
\put(2053,1538){\makebox(0,0)[lb]{{\SetFigFont{12}{14.4}
{\rmdefault}{\mddefault}{\updefault}2}}}
\put(2081,0){\makebox(0,0)[lb]{{\SetFigFont{12}{14.4}
{\rmdefault}{\mddefault}{\updefault}3}}}
\end{picture}
}
\end{minipage}
%%%%%%%
\begin{minipage}{66mm}
\begin{eqnarray*}
A
=
\left[
  \begin{array}{ccccc}
    88 & 50 & 35 & 78 & 38 \\
    50 & 48 & 19 & 27 & 11 \\
    35 & 19 & 41 & 13 & 6 \\
    78 & 27 & 13 & 86 & 44 \\
    38 & 11 & 6 & 44 & 59 \\
  \end{array}
\right]
\end{eqnarray*}
\end{minipage}
\end{center}
%%%%%%%
\end{example}
Here the eigenvector associated with the smallest eigenvalue,
$\lambda_{5}\approx-2.54$, does not have the predicted sign
pattern.
The following is the eigenvector in question, with each entry
approximated to the nearest hundredth:
\begin{eqnarray*}
\textbf{x}
\approx
\left[
  \begin{array}{c}
    -68.08 \\
    32.75 \\
    26.69 \\
    45.57 \\
    \fbox{1} \\
  \end{array}
\right].
\end{eqnarray*}

\section*{Acknowledgements}
We are grateful for the exhaustive comments given by the referee. 
His comments and suggestions have improved the presentation 
of the manuscript. The author R. S. Costas-Santos acknowledges 
financial support by Direcci\'on General de Investigaci\'on, Ministerio 
de Econom\'ia y Competitividad of Spain, grant MTM2012-36732-C03-01.
\bibliographystyle{elsarticle-num}

\begin{thebibliography}{00}

\bibitem{JCB} C. R. Johnson, R. S.
Costas-Santos, and B. Tadchiev.
Matrices Totally Positive Relative to
a Tree. {\it  Electron. J. Linear Algebra} {\bf 18}
(2009), 211--221.

\bibitem{FJ} S. Fallat and C. R. Johnson, {\it
Totally Nonnegative Matrices}, Princeton Series in
Applied Mathematics. Princeton University Press,
Princeton, NJ,
2011. xvi+248 pp.

\bibitem{HoJo13} R. A. Horn and C. R. Johnson, {\it
Matrix Analysis, 2nd Ed.}, Cambridge University Press,
Cambridge, 2013. xviii+643 pp.
\end{thebibliography}

\end{document}